\documentstyle[11pt,amssymb]{article}

\newcommand{\qed}{\relax{\ifhmode\unskip\nobreak\hfill$\Box$\fi
	\ifmmode\ifinner\else\hskip5pt\fi \hfill\Box\fi}\relax}
\newcommand{\Z}{{\Bbb Z}}
\newcommand{\R}{{\Bbb R}}
\newtheorem{theorem}{Theorem}
\newtheorem{question}{Question}
\newenvironment{proof}{{\sc Proof: }}{\qed}

\begin{document}

\title{Another homogeneous, non-bihomogeneous Peano continuum}
\author{Greg Kuperberg}
\date{1997}
\maketitle

\begin{abstract}
K. Kuperberg found a locally connected, finite-dimensional continuum which is
homogeneous but not bihomogeneous. We give a similar but simpler example.
\end{abstract}

In 1988, K. Kuperberg \cite{Mom:knaster} found a locally connected,
finite-dimensional continuum which is homogeneous but not bihomogeneous.
The basic idea of this construction is to use the fact that
a local homeorphism of $\mu^1 \times \mu^1$, where $\mu^1$ is the Menger
curve, either preserves horizontal and vertical fibers or switches
them.  This places strong restrictions on the possible homeomorphisms
of more complicated continua that resemble $\mu^1 \times \mu^1$.
Since $\mu^1$ is homogeneous, the restrictions do not necessarily
preclude homogeneity.

In this note, we give a similar but slightly simpler example of a locally
connected continuum (a Peano continuum) which is homogeneous but not
bihomogeneous.  The author found the example after some informative
discussions with K. Kuperberg.

\section{Terminology}

A space is {\em homogeneous} if for every $x$ and $y$, there is a
homeomorphism sending $x$ to $y$.  It is {\em bihomogeneous} if there is a
homeomorphism exchanging $x$ and $y$.  It is {\em strongly locally
homogeneous} if for every $x$ and every open $U \ni x$, there is
an open $V$ such that $x \in V \subset U$ and for every $y \in V$,
there is a homeomorphism taking $x$ to $y$ which is the identity outside
$U$.  Finally, a space is {\em weakly
2-homogeneous} if for every pair of unordered pairs of distinct points
$\{x_1,x_2\}$ and $\{y_1,y_2\}$, there is homeomorphism taking $\{x_1,x_2\}$
to $\{y_1,y_2\}$.

\section{The example}

Given a closed PL manifold $M$ of dimension at least $2n+1$, let $X_0 = M$,
and for each $i > 0$, let $X_i$ be a regular neighborhood in $X_{i-1}$ of the
$n$-skeleton of a triangulation of $X_{i-1}$. The intersection $\mu^n_M =
\bigcap_i X_i$ is a locally $n-1$-connected continuum called a Menger
manifold whose homeomorphism type does not depend on the particular
triangulations chosen \cite{Bestvina}. The homeomorphism type of $\mu^1_M$
also does not depend on $M$ and is therefore that of $\mu^1$, the Menger
curve \cite{Anderson}. The space $\mu^n_M$ is strongly locally homogeneous.
In general $\pi_k(\mu_M^n) \cong \pi_k(M)$ for $k<n$.  The manifold
$\mu^n_{S^{2n+1}}$ is simply written as $\mu^n$; any $x \in \mu^n_M$ has a
closed neighborhood homeomorphic to $\mu^n$.

Consider the 21-element $G = \Z/3 \ltimes \Z/7$ which is a non-trivial
semidirect product.  The group $G$ is generated by elements $a$ and $b$ with
relations $a^3 = b^7 = 1$ and $aba^{-1} = b^2$. It has no automorphism
sending $a$ to $a^{-1}$.  Let $W$ be a closed, PL 5-manifold with fundamental
group $G$. (One can make a closed, PL 5-manifold whose fundamental group is
an arbitrary finitely presented group by first making a 2-complex from the
presentation, embedding it in $\R^5$, and gluing two copies of a regular
neighborhood along the boundary.) Let $\widetilde{W}$ be the connected covering
of $W$ corresponding to the subgroup generated by $b$, which is a regular
3-fold covering. Then $\pi_1(\widetilde{W}) \cong \Z/7$, $\mu_{\widetilde{W}}^2$ and
$\mu_W^2$ have the same fundamental groups as $\widetilde{W}$ and $W$, and
$\mu_{\widetilde{W}}^2$ is a regular 3-fold covering of $\mu_W^2$.  Evidently,
$\mu_{\widetilde{W}}^2$ has a homeomorphism $\phi$ (a deck translation) of order
three with no fixed points which squares elements of its fundamental group.

The example $X$ is a twisted bundle with fiber $\mu_{\widetilde{W}}^2$ and base
$\mu^1$. Specifically, let $f:\mu^1 \to C$ be a retraction onto a circle $C$
contained in $\mu^1$, let $X'$ be a $\mu_{\widetilde{W}}^2$ bundle over $C$ with
monodromy $\phi$, and let $X$ be the pullback of $X'$ under $f$.
Alternatively, $\mu^1$ has a regular covering of itself of order 3 with deck
translation $\alpha$; let $X$ be the quotient of $\mu^2_{\widetilde{W}} \times
\mu^1$ by the group action $(x,y) \mapsto (\phi(x),\alpha(y))$.

\begin{theorem} The continuum $X$ is homogeneous but not bihomogeneous.
\end{theorem}
\begin{proof}
To establish homogeneity, it suffices to note that $X$ is a regular 3-fold
covering of $\mu^2_W \times \mu^1$.  Let $\pi$ be the covering map. Since
$\mu^2_W$ and $\mu^1$ are strongly locally homogeneous and compact, for every
$x,y \in \mu^2_W \times \mu^1$, there is a homeomorphism of $\mu^2_W \times
\mu^1$ which takes $x$ to $y$ and which lifts to $X$.  Therefore any fiber of
$\pi$ can be taken to any other by a homeomorphism, and it remains to show
that if $x$ and $y$ are in the same fiber, there is a homeomorphism taking
$x$ to $y$.  Indeed, there is a deck translation with this property.

To show that $X$ is not bihomogeneous, we define a horizontal fiber of $X$ to
be the preimage of a point of $f \circ \pi$ and a vertical fiber to be the
preimage of $g \circ \pi$, where $f$ and $g$ are the projections of $\mu^2_W
\times \mu^1$ onto the first and second Cartesian factors, respectively.
By a lemma of K. Kuperberg \cite[lemma 1]{Mom:menger}, arbitrarily small
subsets of $\mu^2 \times \mu^1$ have the property that every open embedding
into must take horizontal fibers (meaning $\mu^2 \times \{p\}$) into
horizontal fibers and vertical fibers (meaning $\{q\} \times \mu^1$) into
vertical fibers.  Therefore a homeomorphism $\psi$ of $X$ must preserve
horizontal and vertical fibers.  In particular, $\psi$ must preserve
$\pi$-fibers, because they are characterized as the intersections of
horizontal and vertical fibers.

We claim that if $x$ and $y$ lie in the same $\pi$-fiber, $\psi$ cannot
exchange them.  Otherwise, $\psi$ fixes the third point of the same fiber.
The map $\psi$ restricts to the vertical fiber containing $x$ and $y$.
Since it preserves $\pi$-fibers, it further descends to a homeomorphism of a
copy of $\mu^2_W$.  Call this homeomorphism $\psi$ also.  Since $\psi$
reverses the cyclic ordering of a $\pi$-fiber over this copy of $\mu^2_W$,
the induced map on $\pi_1(\mu^2_W) \cong G$ sends $a$ to $a^{-1}$.
But $G$ has no such automorphism.
\end{proof}

\section{Remarks and open problems}

Clearly, the 21-element group in the above construction can be replaced by
many other finite groups.  For any $k\ge 1$, there is also a locally
$k$-connected example which is a $\mu^r_W$ bundle over $\mu^s_{S^{2s} \times
S^1}$ for any $r,s > k$.  But if $r=s=1$, the corresponding bundle
becomes bihomogenous.

Say that a topological space is {\em primitively homogeneous} if no
non-trivial equivalence relation on its points is invariant under all
homeomorphisms.  (An equivalence relation is {\em invariant} if that $x
\sim y \iff f(x) \sim f(y)$.  This is analogous to the notion of a
primitive permutation group.)  The continuum $\mu^1 \times \mu^1$ is
primitively homogeneous and bihomogeneous, but not weakly 2-homogeneous.  On
the other hand, if one takes a path-connected component of a solenoid (a
solenoid composant), removes a point, and takes a path-connected component of
the remainder, the result is a topological space embeddable in $\R^n$ which
is primitively homogeneous and even weakly 2-homogeneous, but not bihomogeneous.

\begin{question} Is there a continuum which is primitively homogeneous but not
bihomogeneous?  Primitively homogeneous but neither bihomogeneous nor
weakly 2-homogeneous?
\end{question}


\end{document}